\documentclass{l4dc2024}

\usepackage{comment}
\usepackage{graphicx} % for pdf, bitmapped graphics files
\usepackage{mathptmx} % assumes new font selection scheme installed
\usepackage{times} % assumes new font selection scheme installed
\usepackage{amsmath,bm,mathtools} % assumes amsmath package installed
\usepackage{amssymb}  % assumes amsmath package installed
\usepackage{placeins}
\usepackage{dsfont}
\usepackage{algpseudocode}
\usepackage{algorithm}
\usepackage{url}
\usepackage{cite} %PG : cleans up appearance of multiple citations
\usepackage{calrsfs}
\DeclareMathAlphabet{\pazocal}{OMS}{zplm}{m}{n}
\usepackage{xcolor}
\usepackage{enumitem}

\newtheorem{assumption}{Assumption}
\newtheorem{thrm}{Theorem}
\newtheorem{cor}{Corollary}

\newcommand{\R}{\mathbb{R}}

\newcommand{\bs}[1]{\boldsymbol{#1}}
\newcommand{\zcore}[1]{C_{\zeta}{(#1)}}

\newcommand{\msample}{\boldsymbol{\xi}}

\DeclareMathOperator*{\argmax}{arg\,max}
\DeclareMathOperator*{\argmin}{arg\,min}

\usepackage{xcolor,calc}

\newcommand{\GP}[1]{\textcolor{black}{#1}}
\newcommand{\FF}[1]{\textcolor{black}{#1}}

\title[PAC stability of allocations in coalitional games with private sampling]{Probably approximately correct stability of allocations in uncertain coalitional games with private sampling}

\author{%
	\Name{George Pantazis} \Email{g.pantazis@tudelft.nl}\\
	\addr Delft Center for Systems and Control, TU Delft, The Netherlands
	\AND
	\Name{Filiberto Fele} \Email{ffele@us.es}\\
	\addr Department of Systems Engineering and Automation, University of Seville, Spain%
	\AND
	\Name{Filippo Fabiani} \Email{filippo.fabiani@imtlucca.it}\\
	\addr IMT School for Advanced Studies Lucca, Italy
	\AND
	\Name{Sergio Grammatico} \Email{s.grammatico@tudelft.nl}\\
		\addr Delft Center for Systems and Control, TU Delft, The Netherlands
	\AND
	\Name{Kostas Margellos} \Email{kostas.margellos@eng.ox.ac.uk}\\
	\addr Department of Engineering Science, University of Oxford, UK
}
\begin{document}
\maketitle
\thispagestyle{empty}
\pagestyle{empty}

%%%%%%%%%%%%%%%%%%%%%%%%%%%%%%%%%%%%%%%%%%%%%%%%%%%%%%%%%%%%%%%%%%%%%%%%%%%%%%%%
\begin{abstract} 
We study coalitional games with exogenous uncertainty in the coalition value, in which each agent is allowed to have private samples of the uncertainty. As a consequence, the agents may have a different perception of stability of the grand coalition.
In this context, we propose a novel methodology to study the out-of-sample coalitional rationality of allocations in the set of stable allocations (i.e., the \emph{core}). 
Our analysis builds on the framework of probably approximately correct learning.
Initially, we state \emph{a priori} and \emph{a posteriori} guarantees for the entire core. 
Furthermore, we provide a distributed algorithm to compute a compression set that determines the generalization properties of the a posteriori statements. 
We then refine our probabilistic robustness bounds by specialising the analysis to a single payoff allocation, taking, also in this case, both \emph{a priori} and \emph{a posteriori} approaches.
Finally, we consider a relaxed $\zeta$-core to include nearby allocations and also address the case of empty core. For this case, probabilistic statements are given on the eventual stability of allocations in the $\zeta$-core.
\end{abstract}


\section{Introduction}
Multi-agent systems are pervasive across various fields, including engineering \citep{Raja, KaracaKamgarpour2020TSM, Malcolm2, Fele_coalitional,FiliEtAl2018TAC}, economics, and social sciences \citep{Mccain2008}. Even though agents often behave as self-interested entities, their limited ability to improve their utility can, in certain scenarios, motivate them to form coalitions to achieve a higher individual payoff. This situation can be modelled through coalitional games \citep{Chalkiadakis}.
As each agent's participation in a coalition is subject to their own payoff being maximised, a question emerges about how to allocate the total value of the coalition such that no agent deviates from it.
This problem is known as \emph{stability} of agents' allocations.  \par 
In real-world scenarios, the value of a coalition is typically affected by uncertainty.
The consideration of uncertainty in the coalitional values of a game finds its roots in seminal works such as  \citet{Charnes_Chance_I,Charnes_Chance_II,Charnes1973}, 
which provided solution concepts to yield stable allocations against uncertainty.
\citet{Suijs1} discussed non-emptiness of the core for a particular class of stochastic games. \citet{Chalkiadakis_RL} and \citet{Ieon_Bayesian} addressed uncertainty through Bayesian learning methods, while \citet{Yuqian2015} investigated different solution concepts that maximize the probability of obtaining stable allocations. \citet{Repeated2} and \citet{Repeated1} explored the dynamics of repeated stochastic coalitional games.
 
The connection between probably approximately correct (PAC) learning and uncertain coalitional games has been explored in \citet{Balcan2015}. %In that work, a sampling-based approach is employed to learn the value function, utilizing a randomly generated subset from the total number of potential coalitions.  
The spirit of \citet{Procaccia} is similar, using Vapnik-Chervonenkis (VC) theory to learn the winning coalitions for the class of the so-called \emph{simple} games.  \citet{Balcan2015} focused on a complementary problem where only a randomized subset of coalitions is considered. Both these works evaluate the sample complexity from the VC theory perspective, and hence their results suffer from the associated conservativeness.
%.  As such, these results often suffer from the conservatism associated with VC-theoretic results. 
  \citet{Pantazis_2022_coalitional} leverage the scenario approach  \citep{campi2021theory,campi2018introduction,campi2018waitandjudge,garatti2022risk} to provide  distribution-free guarantees on the stability of allocations in a PAC manner, based on samples from the exogenous uncertainty affecting the value functions. On a parallel line of research, \citet{pantazis2023distributionally} propose a data-driven Wasserstein-based distributionally robust approach for allocations' stability. \par %This implies ensuring, with high confidence, that the probability of agents defecting from the grand coalition and forming subcoalitions is bounded by a predetermined threshold. The uniqueness of our methodology lies in its dependence on the available data, the confidence parameter, and the identification of crucial data samples needed to derive a specific allocation or set of allocations. 
%Notably, the connection between PAC learning and uncertain coalitional games has been explored in \citet{Balcan2015}. In that work, a sampling-based approach is employed to learn the value function, utilizing a randomly generated subset from the total number of potential coalitions. 
%Our work instead shares a similar spirit  with \citet{Pantazis_2022_coalitional} which  leverages recent results from the scenario approach  \citet{Ramponi2018,campi2020scenario,campi2021theory} to provide distribution-free PAC-type stability guarantees for agents'  allocations. Our main contribution compared with \citet{Pantazis_2022_coalitional} is that we extend the former to the more general setting, where agents are allowed to have their own private data sets. 
 Unlike the aforementioned works, we consider a more general setting where the uncertainty data is privately drawn by each agent. As such, information is heterogeneous across agents, and the samples are regarded as a private resource. 
When data samples are commonly shared among agents, the work of \citet{Pantazis_2022_coalitional} provide guarantees on stability of allocations. %, which however cannot be applied to the case where agents draw their own private samples. 
In particular, the developments in \citet{Pantazis_2022_coalitional} rely on the notion of scenario core -- a data-driven approximation of the robust core based on the hypothesis that all uncertainty samples are shared among the agents \citep{Raja}. Given this common set of samples, the perception of a stable allocation based on the available data is identical across agents. This is no longer the case for private sampling, as different samples of the exogenous uncertainty can result in a different perception of stability by each agent. We adopt a PAC learning approach and leverage the concept of compression \citep{MargellosEtAl2015Compression}, i.e., the set of samples essential for the reconstruction of the scenario core and whose cardinality affects the generalization properties of the provided guarantees. We give \emph{a priori} and \emph{a posteriori} certificates for the entire scenario core obtained by private sampling, and propose a distributed algorithm to calculate a compression set.

We then focus on the specific allocation returned by some algorithm (akin, e.g., to the distributed payoff algorithm by \citet{Raja}). We then provide \emph{a priori} coalition stability guarantees for this allocation (Theorem \ref{thm:support_rank_conservative}). %and show that, for the case where there is some underlying structure in the interaction network, i.e., agents are allowed to form subcoalitions only with some but not all agents, one can leverage results on the so called support rank \citep{Schildbach2012} to improve the previous guarantees (see Remark \ref{thrm:support_rank_apriori_improved}). 
We show how the dimension of the problem -- in our case the number of agents -- plays a key role.
% in the provided guarantees. 
 Less conservative probabilistic bounds can often be obtained by taking an \emph{a posteriori} approach \citep{Ramponi2018}, i.e., based on observation of the realized uncertainty. In this case, the probabilistic bounds can be significantly improved and even a tighter \emph{a priori} bound can be obtained (see Theorem \ref{thm:aposteriori_allcoalition}). 
 Finally, we consider a relaxed $\zeta$-core to include nearby allocations and also address the case of empty core. We leverage recents results by \citet{campi2021theory} and study the probabilistic stability of allocations in the $\zeta$-core.

\section{Scenario approach for stochastic coalitional games} \label{formulation}
%\subsection{Problem formulation and stability of allocations}
We consider a coalitional game with $N$ agents identified by the set $\pazocal{N}=\{1, \dots, N\}$.  We denote the number of possible subcoalitions, excluding the grand coalition, by $M$, i.e., $M=2^N-1$. %In this setting, the agents wish to form a certain coalition in case they experience an increase in their individual payoffs. 
The worth of a coalition $S \subseteq \pazocal{N}$ is given by the so called \emph{value function}; here we consider transferable utility problems, i.e., the value of coalitions is expressed by a real value which can be split among participating agents.  We posit that the value of each coalition is subject to uncertainty. % thus rendering the value function  of each coalition stochastic. %Let us, thus define the value function as follows.
Then, the value of a coalition $S \subset \pazocal{N}$ is a function $u_S: \Xi \to \mathbb{R}$ that given an uncertainty realization  $\xi \in \Xi$ returns the total payoff for the agents in $S$.  The value of the grand coalition is the deterministic  quantity $u_\pazocal{N} \in \mathbb{R}$.
%\end{definition}

%\FFnote{Here and in the remainder, the game is defined in several ways, i.e., $G_\mathbb{P}=\langle \pazocal{N},\{u_S\}_{S \subseteq \pazocal{N}}, \Xi, \mathbb{P} \rangle$, or $G_\mathbb{P}=\langle \pazocal{N},\{u_S\}_{S \subseteq \pazocal{N}}, \Xi \rangle$, or $G_\Xi$...}

The uncertain coalitional game is then defined as the tuple $G_\Xi=\langle \pazocal{N},\{u_S\}_{S \subseteq \pazocal{N}}, \Xi, \mathbb{P} \rangle$, where $\mathbb{P}$ is some unknown probability measure over $\Xi$.  A vector  $\bs{x}\coloneqq (x_i)_{i \in \pazocal{N}} \in \mathbb{R}^N$, where $x_i$ is the payoff received by agent $i$, is called an \emph{allocation}.  For a given uncertainty realization $\xi \in \Xi$, an allocation is \emph{strictly rational} for the members of $S$ if $\sum_{i \in S}x_i > u_S(\xi)$; the latter implies that agents have an incentive to form this coalition for this particular uncertainty realization. An allocation is \emph{efficient} if $\sum_{i \in \pazocal{N}}x_i = u_\pazocal{N}$. %In this work we are interested in finding efficient allocations that are not strictly feasible by any other coalition. 
Efficient allocations  such that there are no incentives for agents to deviate the grand coalition are called \emph{stable}. The set of all such allocations is the \emph{core} of the game. 

The notion of robust core is proposed to account for the presence of uncertainty; see, e.g., \citet{Pantazis_2022_coalitional,Raja,NedicBauso2013TAC}. We formally define it as $C(G_\Xi)\coloneqq \{\bs{x} \in  \R^N  \colon \textstyle\sum_ {i \in \pazocal{N}}x_i = u_\pazocal{N},\, \textstyle\sum_{i \in S}x_i \geq \max_{\xi \in \Xi} u_S(\xi) \text{ for all }  S \subset \pazocal{N} \}$.
%$C(G_\Xi)$ provides the coalitional game under study with a measure of robust stability in the sense that, 
For any possible realization of the uncertainty $\xi\in\Xi$, any allocation $\bs{x} \in C(G_\Xi)$ gives the agents no incentive to defect from the grand coalition and form sub-coalitions. Unfortunately, computing explicitly the robust core is hard, as we assume no knowledge on the uncertainty support $\Xi$ (nor on the underlying probability distribution $\mathbb{P}$).
To circumvent this challenge, we adopt a data-driven methodology and approximate the robust core by drawing a finite number $K$ of independent and identically distributed  (i.i.d.) samples $\msample \coloneqq (\xi^{(1)}, \ldots, \xi^{(K)} ) \in \Xi^K$, where $\Xi^K$ denotes the $K$-fold cartesian product of $\Xi$; we refer to vectors $\msample$ as multi-samples. This constitutes the scenario  game $G_K=\langle \pazocal{N}, \{u_S\}_{S \subseteq \pazocal{N}}, \msample \rangle$, whose core is 
	the set $\{ \bs x \in \mathbb{R}^{N}: \sum_{i \in \pazocal{N}} x_i = u_\pazocal{N} \text{ and }\sum_{i \in S}x_i \geq \max_{k=1, \dots K}u_S(\xi^{(k)}),  \, \forall S \subset \pazocal{N} \}$,
referred to as the \emph{scenario core}.

%The first condition is called the efficiency condition. Due to our assumption that the grand coalition is deterministic, it means that the total increase in gains when all agents work together is known with certainty.  
%This is the case when agents might know how efficient a fully-cooperative scheme is but have some level of uncertainty/ambiguity with respect to the potential outcomes of the subcoalitions. The second condition implies that the allocation $x$  is not strictly feasible, hence agents do not have an incentive to form $S$. Otherwise, if $\sum_{i \in S}x_i < u_S(\xi)$  agents would have the incentive to leave the grand coalition and form $S$, thus receiving $ u_S(\xi)$ as opposed to $\sum_{i \in S}x_i$. 

We now take the notion of allocation stability a step further, considering the more general setting where every agent $i$ has only access to a private set of samples $\msample_i$ from $\Xi$. %In this case, even for agents participating in the same coalition, their perception of stability can be different depending on the multi-sample they have drawn. This is formally introduced in the following definition: 
Let $\msample_i \in \Xi^{K_i}$ be the multi-sample privately drawn by agent $i$. We say that an allocation  $\bs x$ of $G_K$ is \emph{stable with respect to} $\msample_i$ if 
$\sum_{i \in \pazocal{N}} x_{i} = u_\pazocal{N}$ and  $\sum_{i \in S}x_{i} \geq \max_{k=1, \dots K_i}u_S(\xi_i^{(k)})$,  $\forall  S \subset \pazocal{N}$  s.t.~$S\supseteq\{i\}$. %\footnote{For brevity, in the following we use $S\supseteq\{i\}$ to denote all $S \subset \pazocal{N}$ which allow agent $i$ as a member.}

%\begin{definition}[Stability w/ private sampling]
%		Let $\msample_i \in \Xi^{K_i}$ be some multi-sample drawn by agent $i$. An allocation  $\bs x=(x_{i})_{i \in \pazocal{N}}$ of the game $G_K$ is stable w.r.t.~$\msample_i$ if 
%$\sum_{i \in \pazocal{N}} x_{i} = u_\pazocal{N}$ and  $\sum_{i \in S}x_{i} \geq \max_{k=1, \dots K_i}u_S(\xi_i^{(k)})$,  $\forall  S \subset \pazocal{N}$  s.t.~$S\supseteq\{i\}$.
%\end{definition}

This immediately leads to the following extension of the scenario core:

\begin{definition} \label{def:scenario_core_private} 
	Let $\msample_i=(\xi_i^{(1)},\ldots,\xi_i^{(K_i)}) \in \Xi^{K_i}$ be some multi-sample drawn by agent $i$. The scenario core with private sampling  is  given by 
	\begin{equation}
		C(G_K)=\bigg\{ \bs x \in \mathbb{R}^{N}: \sum_{i \in \pazocal{N}} x_i = u_\pazocal{N} \text{ and }\sum_{i \in S}x_{i} \geq \max_{i \in S}\max_{k=1, \dots K_i}u_S(\xi^{(k)}),  \, \forall S \subset \pazocal{N} \bigg\}.   \label{eq:scenario_cor_private}
	\end{equation}
\end{definition}

%It is important to highlight 
%There is a crucial conceptual distinction between the standard scenario core and that with private sampling. With the former, all agents share access to the same multi-sample, thereby ensuring a uniform perception of stability among them. Conversely, in the latter scenario, agents lack access to the multi-samples of other agents, as these are treated as private. 
%While coordination is essential to reach a solution that is stable w.r.t. each agents' data, consensus can be achieved in a distributed setting through the corresponding Lagrange multipliers of the coalitional constraints. As such, access to the other agents' samples is not required.%Consequently, if an individual agent  knows of the structure of coalitional value functions with respect to the uncertain parameter, their access to the entire set of multi-samples enables the calculation of all coalitional values appearing in the historical data. \par 

Unless  differently specified, we assume the following conditions hold throughout the paper:
\begin{assumption}\label{ass:iid_nonempty}
	\begin{enumerate}[label=(\roman*)]
		\item \label{ass:iid} Each agent $i \in \pazocal{N}$ draws $K_i$ independent samples from the probability distribution $\mathbb{P}$.  The samples drawn by any agent are independent from those drawn by other agents.
		\item \label{nonempty} $C(G_K)$ is non-empty for any multi-sample $(\msample_i)_{i\in\pazocal{N}} \in \Xi^K$.
	\end{enumerate}
%	i) Each agent $i \in \pazocal{N}$ draws $K_i$ independent samples from the probability distribution $\mathbb{P}$.  The samples drawn by each agent are independent from those drawn by other agents.\\
%	ii) $C(G_K)$ is non-empty for any multi-sample $(\msample_i)_{i\in\pazocal{N}} \in \Xi^K$.
\end{assumption}

%For the subsequent developments we impose the following assumption:
%Moreover, unless differently specified, we will assume the following:
%\begin{assumption} \label{nonempty}
%	$C(G_K)$ is non-empty for any multi-sample $(\msample_i)_{i\in\pazocal{N}} \in \Xi^K$. 
%\end{assumption}

In the following, we let $K=\sum_{i\in \pazocal{N}} K_i$, and $\msample = (\msample_i)_{i=1}^N$. Also, for brevity, we will use $S\supseteq\{i\}$ to denote all $S \subset \pazocal{N}$ which allow agent $i$ as a member.

%Assumption \ref{ass:iid_nonempty}-\ref{nonempty} implies the existence of stable allocations w.r.t. the agents' private samples. 

On the basis of privately available data, we wish to provide guarantees on the probability that allocations $\bs{x} \in C(G_K)$ will remain stable (i.e., within $C(G_K)$) for any future, yet unseen, uncertainty realization. Capitalizing on  \citet{Pantazis2020, Fabiani2020b, Pantazis_2022_coalitional}, we define two probabilistic notions of instability. In particular, the first refers to a particular allocation in the core, whereas the second involves the entire set  $C(G_K)$. 
\begin{definition} \label{def:violation} 
	\begin{enumerate}[label=(\roman*)]
	\item  Let  $V: \R^{N} \rightarrow [0,1]$. For any $\bs{x} \in \mathbb{R}^N$,  $V(\bs{x})\coloneqq \mathbb{P}\{\xi \in \Xi :  \exists  S \subset \pazocal{N},\, \textstyle\sum_{i \in S}x_i <  u_S(\xi)\}$ is the \emph{probability of allocation instability}.
	\item	Let $\mathbb{V}: 2^{\R^{N}} \rightarrow [0,1]$. We call $\mathbb{V}(C(G_K))\coloneqq \mathbb{P}\{\xi \in \Xi :  \exists \bs x \in C(G_K),\, S \subset \pazocal{N}\colon  \textstyle\sum_{i \in S}x_i <  u_S(\xi)\}$ 	\emph{probability of core instability}.   
	\end{enumerate}
\end{definition}

$\mathbb{V}(C(G_K))$ thus denotes the probability with respect to the realizations of $\xi$ that, for some $S$ with value function $u_S(\xi)$,
at least one of the allocations in the scenario core will become unstable, i.e., will be dominated by the option of defecting from the grand coalition to form $S$.  
%By leveraging available data, we then aim at bounding with high confidence 

%\begin{definition}\textup{(Algorithm)} \label{def:algorithm} 
%	A mapping $A: \Xi^K \rightarrow 2^{\mathbb{R}^N}$ that takes as input a multi-sample  $\msample \in \Xi^K$ and returns the scenario core of game $C(G_K)$ is called an \emph{algorithm}. 
%\end{definition} 

To bound the probability of core instability in a PAC fashion, we introduce two key concepts from statistical learning theory, namely the \emph{algorithm} and the \emph{compression set} \citep{MargellosEtAl2015Compression}.
The latter refers to the fact that only a subset of data from $\msample$ may be sufficient to produce the same scenario core. As we will see later this underpins the quality of the provided probabilistic stability guarantees. 
\begin{definition} \label{def:compression} 
	A mapping $\pazocal{A}: \Xi^K \rightarrow 2^{\mathbb{R}^N}$ that takes as input a multi-sample  $\msample \in \Xi^K$ and returns the scenario core of game $C(G_K)$ is called an \emph{algorithm}. 
	With $\mathbb{P}^K$-probability one w.r.t. the choice of $\msample$, a subset  $I \subseteq \{\xi^{(1)}, \dots, \xi^{(K)}\}$  is a \emph{compression set} for $\msample$ if  $\pazocal{A}(I) = \pazocal{A}(\msample)$.%, where $\xi_I$ is a vector whose elements are the samples included in $I$.
\end{definition} 
Any compression set of least cardinality is called \emph{minimal} compression set.  %A procedure that enumerates such a set is called a  \emph{compression function}. 
Another important notion used in our derivations is the \emph{support rank} \citep{Schildbach2012}. % formally defined as follows.
\begin{definition}\label{def:supp_rank}
	Consider the maximal unconstrained subspace $L\in \pazocal{L}$ of a constraint in the form $f(\bs x, \xi) \leq 0$, i.e., $L'\subseteq L$ for all $L'\in \pazocal{L}$, where 
		$\pazocal{L}=\bigcap_{\xi \in \Xi} \bigcap_{\bs x \in \mathbb{R}^{N}}\{ \bar{L} \text{ is a linear subspace in }  \mathbb{R}^N \text{ and } \bar{L} \subset F(\bs x, \xi)\}$,
	with $F(\bs x, \xi)=\{\xi \in \mathbb{R}: f(\bs x+\xi, \delta)=f(\bs x, \xi)\}$. 
	Then, the support rank $\rho$ of this constraint is given by $\rho=N-\text{dim}(L)$.
	\end{definition}
In words, the support rank of an uncertain constraint is equal to the dimension of the problem at hand minus the dimension of the maximal unconstrained space of the constraint.
%The rest of the paper is organized as follows. Section II formulates the problem under study defining the concepts of the expected value core and distributional allocation stability.  Section III provides guarantees that the distributionally robust core is contained with a given confidence within the true expected core of the uncertain coalitional game. In Section IV the asymptotic consistency of the distributionally robust core to the expected value core is shown. Finally, the computational tractability of finding a distributionally stable allocation is discussed, accompanied by numerical simulations. Section V concludes then the paper summarizing our main contributions and proposing directions of future research.
%\emph{Notation}: Denote function $[f]^*: \mathbb{R}^m \rightarrow \mathbb{R}$ as the conjugate function of $f: \mathbb{R}^m \rightarrow \mathbb{R}$, i.e., $[f]^*(y)=\sup_{x \in dom(f)}(y^Tx-f(x))$. Furthermore, function $\| \cdot \|_*: \mathbb{R}^m \rightarrow \mathbb{R}_{\geq}$ denotes the dual norm, while $\sigma_X: \mathbb{R}^m \rightarrow \mathbb{R}$ denotes the conjugate of the characteristic function.

\section{PAC stability guarantees for the scenario core with private sampling}
\subsection{A posteriori collective stability guarantees}
In the following theorem we bound with high confidence the probability that some allocation $\bs x \in C(G_K)$ will become unstable, i.e., that the scenario core -- computed on the basis of $K=\sum_{i\in \pazocal{N}} K_i$ samples -- will be reduced after a new uncertainty realization.
%  Note that in our definition an uncertainty affecting a coalition $S$ can live in a different support set and follow a different probability distribution than the uncertainty affecting  coalition $S' \neq S$. 
%However, the consideration of different probability distributions renders the problem significantly more challenging with respect to studying the stability properties of allocations. 
\begin{thrm} \label{a_posteriori_scenario_core}
	Suppose that each agent \FF{independently draws a} multi-sample	$\bs \xi_{i} \in \Xi^{K_i}$ and \FF{let} $\boldsymbol{\xi} = (\bs \xi_{i})_{i=1}^N$.
	%Denote the cardinality of a compression set as $s_K$,
	\FF{Fix a confidence parameter $\beta\in (0,1)$ and choose $\beta_i>0$ such that $\sum_{i \in \pazocal{N}}\beta_i=\beta$. Consider an algorithm that takes as input $\msample$ and returns the scenario core $C(G_K)$. It holds:}
%	\begin{align}
%	C_p(G_{K+1})=\{ x \in \mathbb{R}^{N}: \sum_{i \in \pazocal{N}} x_i = u_\pazocal{N} \text{ and }\sum_{i \in S}x_i \geq \max\{\max_{k=1, \dots K}u_S(\xi^{(k)}), u_S(\xi)\},  \ \forall \ S \subset \pazocal{N} \}. \nonumber 
%	\end{align}
	\begin{equation}\label{eq:thm1}
	\mathbb{P}^K\Big\{\msample \in \Xi^K: \mathbb{V}(C(G_K))  \leq \FF{\sum_{i \in \pazocal{N}} \epsilon_i(s_{i,K})} \Big\} 	\geq  1- \beta,  
	\end{equation}
	\FF{where $s_{i,K}$ is the cardinality of the subset of a compression relative to agent $i$, quantified a posteriori, and $\epsilon_i$ satisfies}
	\begin{equation}\label{eq:epsilon}
		\epsilon_i(K_i)=1,\quad \sum_{k=1}^{K_i-1}{K_i \choose k}(1-\epsilon_i(k))^{K_i-k}=\beta_i. 
	\end{equation}
\end{thrm} 
\emph{Proof}: \FF{First, note that $\msample\in\Xi^K$ due to Assumption~\ref{ass:iid_nonempty}-\ref{ass:iid}.} Then, the following inequalities hold.
\begin{displaymath}
	\begin{split}
	 & \mathbb{P}^K \bigg\{\msample \in \Xi^K: \mathbb{V}(C(G_K)) \leq \FF{\sum_{i \in \pazocal{N}} \epsilon_i(s_{i,K})} \bigg\} \\
	=& \mathbb{P}^K\bigg\{\msample \in \Xi^K: \mathbb{P}\{ \xi \in \Xi: \exists (i, S \supseteq \{i\}, \bs x): \sum_{i \in S} x_i < u_S(\xi)\} \leq \FF{\sum_{i \in \pazocal{N}} \epsilon_i(s_{i,K})} \bigg\}  \\
	=& \mathbb{P}^K\bigg\{\msample \in \Xi^K: \mathbb{P}\big\{ \bigcup_{i \in \pazocal{N}}\{\xi \in \Xi: \exists (i,S \supseteq \{i\}, \bs x): \sum_{i \in S} x_i < u_S(\xi)\big\} \leq \FF{\sum_{i \in \pazocal{N}} \epsilon_i(s_{i,K})} \bigg\}  \\
	\geq & \mathbb{P}^K\bigg\{\msample \in \Xi^K: \sum_{i \in \pazocal{N}}\mathbb{P}\{\xi \in \Xi:\exists (i,S \supseteq \{i\}, \bs x) : \sum_{i \in S} x_i < u_S(\xi)\} \leq \sum_{i \in \pazocal{N}}\epsilon_i(s_{i,K}) \bigg\} \\
	\geq & \mathbb{P}^K \bigg\{\bigcap_{i \in \pazocal{N}}\{\msample \in \Xi^K: \mathbb{P}\{\xi \in \Xi:  \exists (i,S \supseteq \{i\}, \bs x): \sum_{i \in S} x_i < u_S(\xi)\} \leq \epsilon_i(s_{i,K})\} \bigg\} \nonumber \\
	\geq & 1 - \sum_{i \in \pazocal{N}} \mathbb{P}^K \bigg\{\msample \in \Xi^K: \mathbb{P}\{\xi \in \Xi: \exists (i,S \supseteq \{i\}, \bs x) : \sum_{i \in S} x_i < u_S(\xi)\} > \epsilon_i(s_{i,K}) \bigg\}  
	\end{split} 
\end{displaymath}
%\KM{The first inequality, is this both because of subbaditivity as well as due to the fact that the existence wrt $x$ is dropped? Also how was the sum over $\epsilon(s_{i,K})$ introduced? In the last two inequalities, shouldn't it be $\epsilon_i(s_{i,K})$ rather than $\sum \epsilon_i(s_{i,K})$? Where are $s_{i,K}$ introduced?}
Given some new uncertainty realization $\xi$, let $C(G_{K+\xi})$ designate the scenario core built from the multisample $(\zeta^{(1)},\ldots,\zeta^{(K)},\xi)$. The first equality stems from \FF{Definition~\ref{def:violation}, and expresses the fact that for $C(G_{K+\xi})$ to be a \emph{strict} subset of $C(G_K)$, there must exist some allocation in $C(G_K)$ which, due to $\xi$, violates the rationality condition for some subcoalition $S$, i.e., $\sum_{i \in \pazocal{N}} x_i<u_S(\xi)$, causing the departure of some agent $i \in S$.} The third inequality is obtained by applying the subadditivity property, while the second to last from applying Bonferroni's inequality \citep{KostasEtAl2018TAC,FALSONE202020}. 
Note now that \FF{$C(G_K)$ can be found as the solution of the feasibility problem}
\begin{equation}
	\begin{split}
	&\text{Find all } \bs x \in \mathbb{R}^N   \\
	&\text{ s.t. }\sum_{i \in \pazocal{N}} x_i = u_\pazocal{N}, \\ 
	&\text{ and }\sum_{i \in S}x_i \geq \max_{k=1, \dots K_i}u_S(\xi^{(k)}),  \ \forall \ \FF{S\supseteq\{i\}},\, \forall i \in \pazocal{N},
	\end{split}
\end{equation}
\FF{Let $\mathcal{C}_\xi^i=\{ C_i \subset 2^{R^N}: \sum_{i \in S}x_i \geq u_S(\xi),\, \forall  S\supseteq\{i\}, \forall \bs x \in C_i\}$ be the collection of \FF{(sub)sets of allocations} satisfying the \FF{coalitional} constraints obtained from data corresponding to agent $i \in \pazocal{N}$.
Considering the unique set $\bar{C}_i=\{ \bs x \in \mathbb{R}^N: \sum_{i \in S}x_i \geq \max_{k=1, \dots K_i}u_S(\xi_i^{(k)}),  \, \forall S\supseteq\{i\} \}$, note that $\bar{C}_i \in \bigcap_{k=1}^{K_i} \mathcal{C}_{\xi_i^{(k)}}^i$, i.e., $\bar{C}_i$ satisfies the consistency assumption required by \citet[Th.~1]{Ramponi2018}. From this,}
\begin{multline}
	\mathbb{P}^K\Big\{\msample \in \Xi^K: 	\mathbb{P}\{\xi \in \Xi: \exists (S \supseteq \{i\}, \bs x): \sum_{i \in \pazocal{N}} x_i < u_S(\xi)\} > \epsilon_i(s_{i, K}) \Big\}   \\  
	=\mathbb{P}^K\left\{\msample \in \Xi^K:	\mathbb{P}\{\xi \in \Xi: \bar{C}_i  \notin \mathcal{C}_\xi^i\} > \epsilon_i(s_{i, K})\right\} \leq \beta_i,\nonumber 
\end{multline}
\FF{from which \eqref{eq:thm1} follows.} \hfill $\blacksquare$
%where $\epsilon_i(k): \mathbb{N} \rightarrow [0,1]$ is as defined in \eqref{eq:epsilon}, $\beta_i \in (0,1)$ is a confidence parameter related to each agents' multi-sample which \FF{satisfies $\sum_{i \in \pazocal{N}}\beta_i=\beta$} and \eqref{eq:epsilon}, and $s_{i,K}$ the cardinality of the subset of the compression set \FF{arising from} the randomized coalitional constraints of agent $i$.
%\begin{align}
%	\epsilon(K_i)=1, \sum_{k=1}^{K_i }{N \choose k}(1-\epsilon(k))^{K_i-k}=\beta_i
%\end{align}
%As such, it holds that
%\FF{This concludes the proof, as}
%\begin{displaymath}
%	1 - \sum_{i \in \pazocal{N}} \mathbb{P}^K \{\xi_K \in \Xi^K: \mathbb{P}\{\xi \in \Xi: \exists (S \supseteq \{i\}, x): \sum_{i \in \pazocal{N}} x_i < u_S(\xi)\} \leq \epsilon_i(s_{i,K})\} 
%	\geq  1 - \beta. \qquad \blacksquare 
%\end{displaymath}\par 

\FF{A compression (sub)set $I_i$ originated from agent's $i$ uncertainty samples can be obtained by means of Algorithm 1; its cardinality can then be obtained as $s_{i,K}=|I_i|$.}
\begin{algorithm}[tb]
	\caption{Distributed compression algorithm}\label{alg_distcompr}
	\begin{algorithmic}[1] % This 1 is the starting line number for the algorithm
		\item \textbf{Input:} Multi-sample $\msample_i$, coalition values $\{u_S(\cdot)\}_{S\supseteq\{i\}}$; 
		\item \textbf{Output:} Compression set $I_i$;
		\item \textbf{Initialization:} \FF{$I_i=\varnothing$;}
		\item \textbf{Each agent $i \in \pazocal{N}$ performs}
		\item \quad \textbf{For all \FF{$S' \supseteq \{i\}$}}:
		\begin{equation*}
		\begin{split}
			& \bs x^\ast_{S'} \in \FF{\argmin_{\bs x \geq 0}} \; 0 \\
			\text{s.t. } &\sum_{i \in S'} x_i = \max_{k=1, \dots, K_i}u_{S'}(\xi_i^{(k)}), \\
			& \sum_{i \in S} x_i \geq  \max_{k=1, \dots, K_i}u_S(\xi_i^{(k)}), \forall S \supseteq \{i\} \text{ and } S \neq S'.
		\end{split}
		\end{equation*}
		\item \qquad \textbf{If} $\bs x^\ast_{S'} \neq \varnothing$
		\State \qquad \qquad $I_i \leftarrow I_i \cup \argmax\limits_{k=1, \dots, K_i}u_{S'}(\xi_i^{(k)})$;
		\item \qquad \textbf{End If}
		\item \quad \textbf{End For}
		%\State $s_{i,K}=|I_i|$ %\textcolor{red}{Check that line; what is $s_{K_i}$ and what $S_i$?}
	\end{algorithmic}
\end{algorithm}
Note that differently from the compression algorithm in \citet{Pantazis_2022_coalitional}, Algorithm~\ref{alg_distcompr} is distributed among agents. At each iteration of Algorithm~\ref{alg_distcompr}, a feasibility program is solved where coalitional rationality is enforced with equality for each coalition $S'$ in which agent $i$ is allowed to participate. This allows to identify whether the \FF{sample that maximizes $u_S(\cdot)$ is critical} for the construction of a stable \FF{allocation} set for agent $i \in \pazocal{N}$. \FF{The latter is verified if the problem is feasible, and the sample collected as part of the compression set $I_i$.} Note that unless coordination is imposed among agents, the compression set obtained through Algorithm~\ref{alg_distcompr} is non-minimal.

\subsection{A priori collective stability guarantees}
As a byproduct of our \emph{a posteriori} analysis, an \emph{a priori} bound can be obtained. Specifically one can obtain an a priori bound by considering the worst case value of $\sum_{i \in \pazocal{N}} \epsilon_i(s_{i,K})$ over all possible combinations for which $\sum_{i \in \pazocal{N}}s_{i,K} \leq M$ \FF{(where $M=2^N-1$)}. 
The following theorem provides then an \emph{a priori} bound for the entire core with private \FF{samples.}
\begin{thrm}
	Fix $\beta \in (0,1)$ and choose $\beta_i\in (0,1)$ such that $\sum_{i \in \pazocal{N}}\beta_i=\beta$. \FF{It holds} that
		\begin{equation}
		\mathbb{P}^K\{\msample \in \Xi^K: \mathbb{V}(C(G_K)) \leq \epsilon^\ast \} 	\geq  1- \beta,
		\end{equation}
where \FF{$\epsilon^\ast = \max\{\sum_{i \in \pazocal{N}} \epsilon_i(s_i)\colon \sum_{i \in \pazocal{N}}s_i \leq M,\, s_i\in\mathbb{N}\}$, with $\epsilon_i$ defined as in \eqref{eq:epsilon}, is a quantity independent from the given multi-sample.}
%, obtained as 
%\begin{equation}\label{eq:eps_optimization_program}
%	\begin{split}
%		&\epsilon^\ast = \max_{\{s_i\}_{i \in \pazocal{N}} \in \mathbb{N}^N_{+}} \sum_{i \in \pazocal{N}} \epsilon_i(s_i) \\
%		& \text{subject to } \sum_{i \in \pazocal{N}}s_i \leq M, \nonumber
%	\end{split}
%\end{equation}
\end{thrm} 
\emph{Proof:} %For any multi-sample $\xi_K \in \Xi^K$  it holds that $\sum_{i \in \pazocal{N}}s_i \leq M$. 
\FF{We apply the a posteriori result in  \citet[Th.~1]{Pantazis2020} -- which admits the number of facets of a randomized feasibility domain as an upper bound to the cardinality of the minimal compression set -- to $C(G_K)$, where each facet is associated to some subcoalition in $2^{\pazocal{N}}$. 
Now, let $s_{i,K}$ denote the cardinality of the minimal (sub)compression set relative to agent $i$, quantified \emph{a posteriori}: by \citet[Th.1]{Pantazis2020} it holds that $\sum_{i \in \pazocal{N}}s_{i,K} \leq M$. Then, $\sum_{i \in \pazocal{N}}\epsilon_i(s_{i,K}) \leq \max_{\{s_i\}_{i\in\pazocal{N}}}\sum_{i \in \pazocal{N}} \epsilon_i(s_i)$, for any $\{s_i\}_{i\in\pazocal{N}}$ such that $\sum_{i\in\pazocal{N}} s_i \leq M$. By definition of $\epsilon^*$ it follows}
\begin{equation}
\mathbb{P}^K\{\msample \in \Xi^K:  \mathbb{V}(C(G_K)) \leq \epsilon^\ast \}  \geq \mathbb{P}^K\bigg\{\msample \in \Xi^K:  \mathbb{V}(C(G_K)) \leq \FF{\sum_{i \in \pazocal{N}} \epsilon_i(s_{i,K})} \bigg\} \geq 1 - \beta, \nonumber 
\end{equation}
where the last inequality holds because of Theorem~\ref{a_posteriori_scenario_core}. \hfill $\blacksquare$ \\ 
\FF{The results above can be applied to the entire set of allocations that are stable w.r.t. the observed data; because of their generality, these theoretical guarantees tend to be conservative. In what follows we specialise our analysis to a single allocation.}

\section{PAC stability of \FF{a single} allocation with private sampling }
\FF{Let $\bs x^\ast$ be the unique allocation in $C(G_{K})$ which maximises some convex utility function $f(\cdot)$.\footnote{Uniqueness of the maximiser is without loss of generality. If multiple maximisers are possible, a convex tie-break rule can be applied to single out a unique solution.} Recalling the notion of support rank as per Definition~\ref{def:supp_rank}, the following result holds for $\bs x^*$.}
\begin{thrm} \label{thm:support_rank_conservative}
	Suppose that each agent draws \FF{a multi-sample $\msample_{i}\in\Xi^{K_i}$ and let $\msample = (\msample_i)_{i\in\pazocal{N}}$. Fix $\epsilon \in (0,1)$, and choose $\epsilon_i> 0$ such that $\sum_{i\in\pazocal{N}}\epsilon_i=\epsilon$.
	Then,}
\begin{equation}
	\mathbb{P}^K\big\{\msample \in \Xi^K: \mathbb{P}\{ \xi \in \Xi:  \FF{\bs x^\ast \notin C(G_K)} \} \leq \epsilon \big\} 	\geq  1- \FF{\sum_{i \in \pazocal{N}} \beta_i},  
\end{equation}
\FF{where $\beta_i=\sum\limits_{j=1}^{\rho_i}{K_i \choose j}\epsilon^j_i(1-\epsilon_i)^{K_i-j}$,
with $\rho_i\leq N$ being the support rank of the coalitional rationality constraints corresponding to agent $i$, i.e., relative to all $S\subset \pazocal{N}$ such that $S\supseteq\{i\}$.}
\end{thrm}
\emph{Proof}: 
\FF{Similarly} to the proofline of Theorem 1, we have that 
\begin{align}
	& \mathbb{P}^K\Big\{\msample \in \Xi^K:   V(\bs x^\ast) \leq \epsilon \Big\} \nonumber \\
	=& \mathbb{P}^K\Big\{\msample \in \Xi^K: \mathbb{P}\{ \xi \in \Xi: \exists  (i \in \pazocal{N}, S \supset \{i\}) : \sum_{i \in S} x^\ast_{i} < u_S(\xi)\} \leq \FF{\sum_{i\in\pazocal{N}}\epsilon_i} \Big\} \nonumber \\
	=& \mathbb{P}^K\Big\{\msample \in \Xi^K: \mathbb{P}\{ \bigcup_{i \in \pazocal{N}}\{\xi \in \Xi: \exists S \supseteq \{i\}: \sum_{i \in S} x^\ast_{i} < u_S(\xi)\} \}\leq \FF{\sum_{i\in\pazocal{N}}\epsilon_i} \Big\} \nonumber \\
	\geq & \mathbb{P}^K\Big\{\msample \in \Xi^K: \sum_{i \in \pazocal{N}}\mathbb{P}\{\xi \in \Xi: \exists S \supseteq \{i\}: \sum_{i \in S}x^\ast_{i} < u_S(\xi)\} \leq \FF{\sum_{i\in\pazocal{N}}\epsilon_i} \Big\} \nonumber \\
	\geq & \mathbb{P}^K \Big\{\bigcap_{i \in \pazocal{N}}\{\msample \in \Xi^K: \mathbb{P}\{\xi \in \Xi: \exists S \supseteq \{i\}: \sum_{i \in S} x^\ast_{i} < u_S(\xi)\} \leq \epsilon_i \} \Big\} \nonumber \\
	\geq & 1 - \sum_{i \in \pazocal{N}} \mathbb{P}^K \Big\{\msample \in \Xi^K: \mathbb{P}\{\xi \in \Xi: \exists S \supseteq \{i\}: \sum_{i \in S} x^\ast_{i} < u_S(\xi)\} > \epsilon_i\Big\}
	\geq  1 - \sum_{i \in \pazocal{N}} \beta_i. \nonumber 
\end{align}
To obtain the last inequality consider the following optimization program 
\begin{equation}\label{eq:opt_allocation}
	\begin{split}
	& \FF{\bs x^* = \argmin_{\bs x \in \mathbb{R}^{N}} f(\bs x)}   \\
	&\text{ s.t. }\sum_{i \in \pazocal{N}} x_i = u_\pazocal{N},  \\ 
	&\text{ and }\sum_{i \in S}x_i \geq \max_{k=1, \dots K_i}u_S(\xi_i^{(k)}),\;  \forall \FF{S\supseteq\{i\},}\; \forall i \in \pazocal{N}.
	\end{split}
\end{equation}
 We then group the constraints depending on which agent they correspond to. \FF{Let $X_i(\xi)\coloneqq\{\bs x \in \mathbb{R}^N\colon \sum_{i \in S}x_i \geq u_S(\xi),\, \forall S \supseteq \{i\} \}$.
With $\beta_i$ defined as above, from \citet{Schildbach2012} we have
\begin{equation*}
	\begin{split}
	\beta_i  &{}\geq   \mathbb{P}^K\{\msample \in \Xi^K:    \mathbb{P}\{ \xi \in \Xi: \bs x^\ast  \notin	X_i(\xi)\} > \epsilon_i \}  \\
	&{}=  \mathbb{P}^K\{\msample \in \Xi^K: \mathbb{P}\{ \xi \in \Xi: \exists  (i \in \pazocal{N}, S \supseteq \{i\}) : \sum_{i \in S} x^\ast_{i} < u_S(\xi)\} > \epsilon_i\},   
	\end{split}
\end{equation*}
which concludes the proof.} \hfill $\blacksquare$

Note that the confidence parameter $\beta_i$ in Theorem \ref{thm:support_rank_conservative} depends on the number of samples $K_i$ of agent $i \in \pazocal{N}$ and the violation level $\epsilon_i$ set individually by the agent to determine the probability of satisfaction of the rationality constraints. Most importantly, $\beta_i$ depends on the so called support rank $\rho_i$ of the coalitional constraints corresponding to $i \in \pazocal{N}$, \FF{which is in any instance upper bounded by the dimension of the decision variable in \eqref{eq:opt_allocation}.} %, that is the dimension of the space of the constraints $X^{(i)}_\xi$ minus the dimension of the maximal unconstrained subspace. 
Theorem \ref{thm:support_rank_conservative}  is an  \emph{a priori} result and in certain cases can be conservative. The following result provides an improved probabilistic bound, \FF{based on the observation of the uncertainty realization.}
Let
\begin{equation}\label{eq:epsilon_beta_i}
	\epsilon_i(s)=1 - \Bigg(\frac{\beta_i}{(N+1) {K_i \choose s}}\Bigg)^{\frac{1}{K_i-s}}. 
\end{equation}
Note that \eqref{eq:epsilon_beta_i} is compatible with \eqref{eq:epsilon}, and allows to explicitly compute $\epsilon_i$ as a function of $\beta_i$.
\begin{thrm} \label{thm:aposteriori_allcoalition}
	Fix $\beta \in (0,1)$ and choose $\beta_i>0$ such that $\sum_{i \in \pazocal{N}} \beta_i = \beta$. Set $\epsilon_i$ as in \eqref{eq:epsilon_beta_i}. Then, for the unique allocation $\bs x^*\in C(G_K)$ it holds
		\begin{equation}
		\mathbb{P}^K\bigg\{\msample \in \Xi^K: \mathbb{P}\{ \xi \in \Xi:  \FF{\bs x^\ast \notin C(G_K)} \} \leq \sum_{i \in S}\epsilon_i(s_{i, K}) \bigg\} 	\geq  1- \beta.
	\end{equation}
\end{thrm}
Recall that, as in Theorem~\ref{thm:aposteriori_allcoalition}, $s_{i,K}$ is the cardinality of the (sub)compression set associated with the rationality constraints relative to all $S\supseteq\{i\}$, quantified a posteriori, e.g., by means of the procedure illustrated in \citet[\S II]{Ramponi2018}.

\noindent\emph{Proof}: For each agent $i \in \pazocal{N}$, taking the conditional probability over \FF{all other agents' multi-samples, Theorem 1 in \citet{Ramponi2018} allows to state}
\begin{equation*}
	\mathbb{P}^K\big\{\msample \in \Xi^K: \mathbb{P}\{ \xi \in \Xi:  \sum_{i \in S}x_{i}^\ast > u_S(\xi) \} \leq \epsilon_i(s_{i, K}) |  (\msample_j)_{j \neq i} \in \Xi^{K-K_i}  \big\} 	\geq  1- \beta_i,
\end{equation*}
which, by integrating w.r.t.~the probability of \FF{realization} of all other agents' scenarios, becomes
\begin{align}
	\mathbb{P}^K\big\{\msample \in \Xi^K: \mathbb{P}\{ \xi \in \Xi:  \sum_{i \in S}x_{i}^\ast > u_S(\xi) \} \leq \epsilon_i(s_{i, K}) \big\} 	\geq  1- \beta_i. \label{relation_aposteriori}
\end{align}
\FF{The proof is concluded by} applying  relation (\ref{relation_aposteriori}) in the proofline of Theorem \ref{thm:support_rank_conservative}. \hfill $\blacksquare$

\FF{Finally, an a priori bound can be derived from Theorem~\ref{thm:aposteriori_allcoalition} by considering the worst case among all a posteriori observable compressions $\{s_{i, K}\}$; this bound is nontrivial by noticing that $\sum_{i \in \pazocal{N}}s_{i, K} \leq N$, which follows from \citet[Th.~3]{Calafiore2006}.
\begin{cor}
	Let $\epsilon^\ast = \max\{\sum_{i \in \pazocal{N}} \epsilon_i(s_i)\colon \sum_{i \in \pazocal{N}}s_i \leq N,\, s_i\in\mathbb{N}\}$, with $\epsilon_i$ defined as in \eqref{eq:epsilon_beta_i}. Then
		\begin{equation}
		\mathbb{P}^K\big\{\msample \in \Xi^K: \mathbb{P}\{ \xi \in \Xi:  \FF{\bs x^\ast \notin C(G_K)} \} \leq \epsilon^\ast \big\} 	\geq  1- \beta.
	\end{equation}
\end{cor}}
%where $\epsilon_i(s_{i, K})$ satisfies constraint \eqref{epsilon_constraint} \FF{and $\sum_{i \in \pazocal{N}}s_{i, K} \leq N$ is due to \citet[Thm.~3]{Calafiore2006}.} \hfill $\blacksquare$

\subsection{The case of empty core}
%There might exist realizations of the uncertainty $\xi \in \Xi$ for which we obtain an empty (scenario) core. 
%the same might happen even for the case of the scenario core $C_p(G_K)$. 
%Establishing whether a private scenario core is empty amounts to solving the following optimization problem:
%\begin{equation}
%	\begin{split}
%	& \min_{\bs x \in \mathbb{R}^{N}} 0  \nonumber \\
%	&\text{ s.t. }\sum_{i \in \pazocal{N}} x_i = u_\pazocal{N} \nonumber \\ 
%	&\text{ and }\sum_{i \in S}x_i \geq \max_{k=1, \dots K_i}u_S(\xi_i^{(k)}),  \; \forall S\supseteq\{i\},\; \forall i \in \pazocal{N}
%	\end{split}
%\end{equation}
Lifting Assumption \ref{ass:iid_nonempty}-\ref{nonempty} on non-emptiness of the scenario core, we define a relaxed version of the latter, the so called $\zeta$-core for the case of private uncertainty sampling as follows:
%\KM{In all this section I think $k_i$ should be dropped, and we should have $\zeta_i$ rather than $\zeta$; changes aren't in red but they appear in multiple occasions}
\begin{definition} \label{def:epsilon_scenario_core}
	\FF{For any $(\msample_i)_{i\in\pazocal{N}} \in \Xi^K$, fixed $\bar{\zeta}_i \geq 0$ for all $i \in \pazocal{N}$, the scenario $\zeta$-core is given by}
		\begin{equation}
		\zcore{G_K}=\Big\{ \bs x \in \mathbb{R}^{N}: \sum_{i \in \pazocal{N}} x_i = u_\pazocal{N} \text{ and }\sum_{i \in S}x_{i} \geq \max_{i \in S}\max_{k=1, \dots, K_i}u_S(\xi_i^{(k)})-\bar{\zeta}_i,  \ \forall S \subset \pazocal{N} \Big\}.   \label{eq:zeta_cor_private}
	\end{equation}
\end{definition}
\FF{It is worth pointing out that the $\zeta$-core allows to extend the analysis to allocations ``closest'' to the core, and the interest in it may not be necessarily restricted to cases where the standard core is empty (this is recovered by setting $\bar{\zeta}_i=0$ for all $i \in \pazocal{N}$). Also, Definition~\ref{def:epsilon_scenario_core} contemplates a different relaxation parameter $\bar{\zeta}_i$ for every agent, allowing to specialise the definition according to individual preferences, or possibly different information on, e.g., each agent's sampling.}
%Definition \ref{def:epsilon_scenario_core} allows to use a different parameter $\zeta_i^{(k)}$ for each agent's sample,
%yielding different penalizations depending on how much a particular sample prevents the agents from forming a stable allocation. 

\FF{On these grounds, here we address the problem of i) providing an allocation with (relaxed) coalitional stability certificates and ii) measuring how agents' multi-samples contribute to lack of coalitional} stability. These questions can be answered at once by solving
\begin{subequations}\label{eq:prob_zeta_core}
\begin{align}
	& \min_{\bs x, \{\zeta_i\geq 0\}_{i\in\pazocal{N}}}\; \sum_{i \in \pazocal{N}} \sum_{k=1}^{K_i} \FF{\zeta_i^{(k)} }, \\
	& \text{s.t. } \sum_{i \in \pazocal{N}} x_i = u_\pazocal{N}  \\
	& \text{and } \sum_{i \in S} x_i \geq \max_{k=1, \dots, K_i} u_S(\xi^{(k)}) - \zeta_i^{(k)}, \; \forall  S \supseteq \{i\},\; \forall i \in \pazocal{N}. \label{eq:prob_zeta_core_3} 
\end{align}
\end{subequations}
%Problem \eqref{eq:prob_zeta_core} can be solved in a distributed manner both w.r.t. $\bs x = (x_i)_{i\in\pazocal{N}}$ as well as $\{\zeta_i\}_{i\in\pazocal{N}}$; consensus is needed to satisfy the equality constraint imposed by the grand coalition.

In what follows, we consider without loss of generality that for any multi-sample $(\msample_i)_{i\in\pazocal{N}} \in \Xi^K$, \eqref{eq:prob_zeta_core} returns a unique pair $(\bs x^\ast, \bs \zeta^\ast)\in\mathbb{R}^N\times\mathbb{R}^K$ (this can be ensured by using a convex tie-break rule). 
\begin{assumption} \label{ass:non_accumulation}
	For every allocation $\bs x \in \mathbb{R}^N$, $\mathbb{P}\{ \xi \in \Xi: \sum_{i \in S}x_i=u_S(\xi) \}=0$ for any $S \subset \pazocal{N}$. 
	\end{assumption}
Assumption~\ref{ass:non_accumulation} is related to non-degeneracy and is often satisfied when $\xi$ does not accumulate as is the case for continuous probability distributions.
Now, for $i \in \pazocal{N}$, consider the polynomial equation \citep[Th. 1]{campi2021theory}
\begin{equation}\label{eq1}
	\binom{K_i}{s}t^{K_i-s}-\frac{\beta_i}{2N}\sum\limits_{j=s}^{K_i-1}\binom{i}{s}t^{j-s}-\frac{\beta_i}{6K_i}\sum\limits_{j=K_i+1}^{4K_i}\binom{j}{s}t^{j-s}=0
\end{equation}	
%\begin{cases}
%	\binom{K_i}{s}t^{K_i-s}-\frac{\beta_i}{2N}\sum\limits_{j=s}^{K_i-1}\binom{i}{s}t^{j-s}-\frac{\beta_i}{6K_i}\sum\limits_{j=K_i+1}^{4K_i}\binom{j}{s}t^{j-s}=0, & \text{if $s\in\{0, \dots, K_i-1\}$},   \\
%	1-\frac{\beta_i}{6K_i}\sum\limits_{j=K_i+1}^{4K_i}\binom{j}{K_i}t^{j-K_i}=0, & \text{if } s=K_i. 
%\end{cases}
For $s\in\{0, \dots, K_i-1\}$, \eqref{eq1} has two solutions in $[0, +\infty)$. We denote the smallest as $\underline{t}_i(s)$, and further let $\underline{t}_i(s)=0$ for $s=K_i$. %Then, we define  $\overline{\epsilon}_i(s)=1-\underline{t}_i(s)$, $s=0,1, \dots, K_i$.
We can now propose the following a posteriori statement on $\bs x^*$:
\begin{thrm}\label{thm:zeta_alloc}
	Fix $\beta \in (0,1)$ and choose $\beta_i >0$  such that $\sum_{i \in \pazocal{N}}\beta_i = \beta$. Let $(\bs x^*, \bs \zeta^*)$ be the solution of \eqref{eq:prob_zeta_core}, and denote by $s_i^*$ the cardinality of $\{k\colon \zeta_i^{*(k)}>0\}$. Under Assumption \ref{ass:non_accumulation},
	\begin{align}
		\mathbb{P}^K\big\{\msample \in \Xi^K&: \mathbb{P}\{ \xi \in \Xi:  \exists (i, S \supseteq \{i\}) : \sum_{i \in S}  x^\ast_{i} < u_S(\xi) \} \leq \sum_{i \in \pazocal{N}}\bar{\epsilon}_i(s_i^*)
	\big\} \geq  1- \beta \label{eq:zeta_core},
\end{align}
where $\overline{\epsilon}_i(s_i^*)=1-\underline{t}_i(s_i^*)$, for all $i\in\pazocal{N}$.
	\end{thrm}
\emph{Proof}: Let $f_i(\bs x, \xi)\coloneqq \max_{j = 1,\ldots,M}(b_j(\xi)-a^\top_j \bs x)$, where the $i$-th component of $a_{j}$ is 1 if $i\in S$ (and 0 otherwise), and $b_j(\xi)=(u_S(\xi))_{S \supseteq \{i\}}$. Then \eqref{eq:prob_zeta_core_3} can be rewritten as %$f_i(\bs x, \xi_i^{(k)}) \leq \zeta_i^{(k)}$
\begin{equation*}
	f_i(\bs x, \xi_i^{(k)}) \leq \zeta_i,  \;  \forall k=1,\dots, K_i, \; \forall i \in \pazocal{N}.
\end{equation*}
We note that in our setting \citet[Assum.~1]{campi2021theory} is satisfied; under Assumption~\ref{ass:non_accumulation}, it is then possible to apply \citet[Th.~1]{campi2021theory}, which yields 
	\begin{equation}
	\mathbb{P}^K \left\{ \msample \in \Xi^K: \mathbb{P}\{ \xi \in \Xi: f_i(\bs x^*, \xi) > 0 \} >  \bar{\epsilon}_i(s_i^*) \right\} \leq  \beta_i,
\end{equation}
where $\bar{\epsilon}_i(\cdot)$ is derived from \eqref{eq1} as described above. By applying arguments similar to those used in proving Theorems \ref{thm:support_rank_conservative} and \ref{thm:aposteriori_allcoalition}, we obtain (\ref{eq:zeta_core}). \hfill $\blacksquare$

The statement in Theorem~\ref{thm:zeta_alloc} can be interpreted as follows: as a result of \eqref{eq:prob_zeta_core}, $\bs x^*$ is an allocation in the $\zeta$-core defined by $\{\bar{\zeta}_i\}_{i\in\pazocal{N}}$, with $\bar{\zeta_i}=\max_{k=1,\ldots,K_i} \zeta_i^{*(k)}$. It holds with confidence $1-\beta$ that $\bs x^*$ is a coalitionally stable allocation with probability $1-\sum_{i\in\pazocal{N}}\bar{\epsilon}_i$, under all possible realizations of the uncertain parameter $\xi\in\Xi$.

\section{Conclusion}
In this work we considered uncertain coalitional games and proposed a data-driven methodology to study the stability of allocations for the general setting where uncertainty is privately sampled by the agents. Future work will investigate stability of allocations through a distributionally robust framework, as well as analyse the case where the value of the grand coalition can also be affected by uncertainty. 
	
\acks{F.~Fele gratefully acknowledges support from grant RYC2021-033960-I funded by MCIN/AEI/ 10.13039/501100011033 and European Union NextGenerationEU/PRTR, as well as from grant PID2022-142946NA-I00 funded by MCIN/AEI/ 10.13039/501100011033 and by ERDF A way of making Europe.}

\bibliography{biblio_coop}
\end{document}